April 21, 2013

# EVALUATION OF ONE EXOTIC FURDUI TYPE SERIES

**Khristo N. Boyadzhiev**
Department of Mathematics and Statistics, Ohio Northern University, Ada, OH 45810
k-boyadzhiev@onu.edu

**Abstract**. In this note we present a closed form evaluation of one interesting alternating series whose value is a combination of three mathematical constants, $\pi$, $\ln 2$, and $\zeta(3)$. The series is expressed as a symmetrical double integral and this integral is evaluated in several steps.

## 1. Introduction

The following two series were evaluated recently

$$\sigma_1 = \sum_{n=0}^{\infty}(-1)^n \left(\sum_{k=1}^{\infty}\frac{(-1)^{k-1}}{n+k}\right)^2 = \frac{\pi^2}{24},$$

(problem 11682 in the American Mathematical Monthly [4]) and the more challenging one

$$\sigma_2 = \sum_{n=1}^{\infty}(-1)^n \left(\sum_{k=1}^{\infty}\frac{(-1)^{k-1}}{2n+k}\right)^2 = \frac{G}{2} + \frac{\pi^2}{48} - \frac{7}{8}(\ln 2)^2 - \frac{\pi}{8}\ln 2,$$

where $G$ is the Catalan constant. The evaluation of $\sigma_2$ was published in [2]. This series is remarkable with its connection to the three important constants $\pi$, $\ln 2$, and $G$.

In this note we shall evaluate another interesting series

$$\sigma = \sum_{n=1}^{\infty}\frac{(-1)^{n-1}}{n}\left(\sum_{k=1}^{\infty}\frac{(-1)^{k-1}}{n+k}\right)^2 = \sum_{n=1}^{\infty}\frac{(-1)^{n-1}}{n}\left(\ln 2 - \sum_{k=1}^{n}\frac{(-1)^{k-1}}{k}\right)^2.$$

Namely, we have the following.

**Proposition.** The series $\sigma$ converges and its value is

$$\sigma = \frac{\pi^2}{12}\ln 2 + \frac{1}{3}(\ln 2)^3 - \frac{1}{2}\zeta(3) \approx 0.08007,$$

involving the three constants $\pi$, $\ln 2$, and $\zeta(3)$.

The above series can be expressed as double integrals. Thus

$$\int_0^1\int_0^1 \frac{dx\,dy}{(1+xy)(1+x)(1+y)} = \sigma_1 = \frac{\pi^2}{24},$$

$$\int_0^1\int_0^1 \frac{x^2 y^2\,dx\,dy}{(1+x^2 y^2)(1+x)(1+y)} = -\sigma_2 = \frac{7}{8}(\ln 2)^2 + \frac{\pi}{8}\ln 2 - \frac{G}{2} - \frac{\pi^2}{48},$$

$$\int_0^1\int_0^1 \frac{\ln(1+xy)}{(1+x)(1+y)}\,dx\,dy = \sigma = \frac{\pi^2}{12}\ln 2 + \frac{1}{3}(\ln 2)^3 - \frac{1}{2}\zeta(3).$$

For related results see [3].

## 2. Evaluation of the series

Using the representation

$$\frac{1}{n+k} = \int_0^1 x^{n+k-1}\,dx,$$

we write

$$\sum_{k=1}^\infty \frac{(-1)^{k-1}}{n+k} = \int_0^1 x^n \left\{\sum_{k=1}^\infty (-1)^{k-1} x^{k-1}\right\} dx = \int_0^1 \frac{x^n}{1+x}\,dx.$$

Then

$$\sigma = \sum_{n=1}^\infty \frac{(-1)^{n-1}}{n}\left(\int_0^1 \frac{x^n\,dx}{1+x}\right)^2 = \sum_{n=1}^\infty \frac{(-1)^{n-1}}{n}\left(\int_0^1 \frac{x^n\,dx}{1+x}\right)\left(\int_0^1 \frac{y^n\,dy}{1+y}\right) = \sum_{n=1}^\infty \frac{(-1)^{n-1}}{n}\int_0^1\int_0^1 \frac{x^n y^n\,dx\,dy}{(1+x)(1+y)}$$

$$= \int_0^1\int_0^1 \left\{\sum_{n=1}^\infty \frac{(-1)^{n-1}(xy)^n}{n}\right\} \frac{dx\,dy}{(1+x)(1+y)} = \int_0^1\int_0^1 \frac{\ln(1+xy)}{(1+x)(1+y)}\,dx\,dy$$

$$= \int_0^1 \frac{1}{1+x}\left\{\int_0^1 \frac{\ln(1+xy)}{1+y}\,dy\right\} dx.$$

Now set

$$H_n^- = 1 - \frac{1}{2} + \frac{1}{3} + \ldots + \frac{(-1)^{n-1}}{n}.$$

We shall evaluate the inside integral by using series expansion. For $|x| < 1$ we write

$$\int_0^1 \frac{\ln(1+xy)}{1+y} dy = \sum_{n=1}^{\infty} \frac{(-1)^{n-1} x^n}{n} \left\{ \int_0^1 \frac{y^n}{1+y} dy \right\}$$

$$= \sum_{n=1}^{\infty} \frac{(-1)^{n-1} x^n}{n} \left\{ \sum_{k=1}^{\infty} \frac{(-1)^{k-1}}{n+k} \right\} = \sum_{n=1}^{\infty} \frac{(-1)^{n-1} x^n}{n} \left( \ln 2 - H_n^- \right)(-1)^n$$

and so,

$$\int_0^1 \frac{\ln(1+xy)}{1+y} dy = \ln 2 \ln(1-x) + \sum_{n=1}^{\infty} \frac{x^n}{n} H_n^-. \qquad (1)$$

***Remark 1.*** At this point it is good to notice that a slight modification of entry (5.5.27) in [5] yields

$$\ln 2 \ln(1-x) + \sum_{n=1}^{\infty} \frac{x^n}{n} H_n^- = \text{Li}_2\left(\frac{1-x}{2}\right) - \text{Li}_2(-x) + \frac{1}{2}(\ln 2)^2 - \frac{\pi^2}{12},$$

and as a result from (1) we have the closed form evaluation

$$\int_0^1 \frac{\ln(1+xy)}{1+y} dy = \text{Li}_2\left(\frac{1-x}{2}\right) - \text{Li}_2(-x) + \frac{1}{2}(\ln 2)^2 - \frac{\pi^2}{12} \qquad (2)$$

true on the entire unit disk $|x| \leq 1$. Integrals of this type have been considered in [1] and in Levin's book [6, A.3.1 (2)], but the form given here is simpler and more appropriate for our work.

Now we can write

$$\sigma = \int_0^1 \frac{1}{1+x} \left\{ \int_0^1 \frac{\ln(1+xy)}{1+y} dy \right\} dx$$

$$= \int_0^1 \frac{1}{1+x} \text{Li}_2\left(\frac{1-x}{2}\right) dx - \int_0^1 \frac{1}{1+x} \text{Li}_2(-x) dx + \left\{ \frac{1}{2}(\ln 2)^2 - \frac{\pi^2}{12} \right\} \int_0^1 \frac{dx}{1+x}$$

and we shall evaluate these integrals one by one. Making first the substitution $t = \frac{1}{2}(1-x)$ and then integrating by parts,

$$J_1 = \int_0^1 \frac{1}{1+x} \operatorname{Li}_2\left(\frac{1-x}{2}\right) dx = \int_0^{1/2} \frac{\operatorname{Li}_2(t)}{1-t} dt = \ln 2 \operatorname{Li}_2\left(\frac{1}{2}\right) - \int_0^{1/2} \frac{\ln^2(1-t)}{t} dt.$$

We now use the antiderivative

$$\int_0^t \frac{\ln^2(1-t)}{t} dt = \ln(t)\ln^2(1-t) + 2\ln(1-t)\operatorname{Li}_2(1-t) - 2\operatorname{Li}_3(1-t) + 2\zeta(3) \qquad (3)$$

(see [5, (6.27), p.159]) and setting $t = 1/2$ we find

$$\int_0^{1/2} \frac{\ln^2(1-t)}{t} dt = -(\ln 2)^3 - 2\ln 2 \operatorname{Li}_2\left(\frac{1}{2}\right) - 2\operatorname{Li}_3\left(\frac{1}{2}\right) + 2\zeta(3).$$

It is well known that

$$\operatorname{Li}_2\left(\frac{1}{2}\right) = \frac{\pi^2}{12} - \frac{(\log 2)^2}{2},$$

$$\operatorname{Li}_3\left(\frac{1}{2}\right) = \frac{7}{8}\zeta(3) - \frac{\pi^2 \log 2}{12} + \frac{(\log 2)^3}{6},$$

and we compute

$$J_1 = \frac{\pi^2}{12}\ln 2 - \frac{1}{6}(\ln 2)^3 - \frac{1}{4}\zeta(3).$$

Next we continue with

$$J_2 = -\int_0^1 \frac{\operatorname{Li}_2(-x)\, dx}{1+x} = \int_0^{-1} \frac{\operatorname{Li}_2(t)}{1-t} dt = -\ln 2 \operatorname{Li}_2(-1) - \int_0^{-1} \frac{\ln^2(1-t)}{t} dt.$$

Here $\operatorname{Li}_2(-1) = \frac{-\pi^2}{12}$ is easy, but for the integral on the RHS we cannot use the antiderivative in (3), because it brings to the expressions $\ln(-1), \operatorname{Li}_2(2), \operatorname{Li}_3(2)$ which need to be defined additionally. Instead of discussing these values we shall follow a different approach. Let

$$H_n^{(2)} = 1 + \frac{1}{2^2} + \frac{1}{3^2} + \ldots + \frac{1}{n^2}.$$

Expanding $(1-t)^{-1}$ in power series and using Cauchy's rule for the product of two power series we obtain

$$\frac{\operatorname{Li}_2(t)}{1-t} = \sum_{n=1}^{\infty} H_n^{(2)} t^n,$$

and therefore,

$$J_2 = \int_0^{-1} \frac{\operatorname{Li}_2(t)}{1-t}\, dt = \sum_{n=1}^{\infty} \frac{(-1)^{n+1}}{n+1} H_n^{(2)} = \frac{\pi^2}{12}\ln 2 - \frac{1}{4}\zeta(3).$$

The evaluation of this series is given in the next section.

Finally we have

$$\sigma = J_1 + J_2 + \left(\frac{1}{2}(\ln 2)^2 - \frac{\pi^2}{12}\right)\ln 2$$

and after simplifying we come to the desired result

$$\sigma = \frac{\pi^2}{12}\ln 2 + \frac{1}{3}(\ln 2)^3 - \frac{1}{2}\zeta(3).$$

*Remark 2.* As byproducts from (1) and (2) with $x = -1$ we have

$$\int_0^1 \frac{\ln(1-x)}{1+x}\, dx = \frac{1}{2}(\ln 2)^2 - \frac{\pi^2}{12},$$

$$\sum_{n=1}^{\infty} \frac{(-1)^{n-1}}{n} H_n^{-} = \frac{1}{2}(\ln 2)^2 + \frac{\pi^2}{12}.$$

## 3. Additional results

*Lemma.* For $|t| \leq 1$, $t \neq 1$, we have

$$\sum_{n=1}^{\infty} H_n^{(2)} \frac{t^n}{n} + 2\sum_{n=1}^{\infty} H_n \frac{t^n}{n^2} = 3\operatorname{Li}_3(t) - \operatorname{Li}_2(t)\ln(1-t) .\qquad(4)$$

*Proof.* We already know that

$$\sum_{n=1}^{\infty} H_n^{(2)} t^n = \frac{\operatorname{Li}_2(t)}{1-t}.$$

From this, dividing by $t$ and integrating,

$$\sum_{n=1}^{\infty} H_n^{(2)} \frac{t^n}{n} = \int_0^t \frac{\text{Li}_2(t)}{t(1-t)} dt = \int_0^t \frac{\text{Li}_2(t)}{t} dt + \int_0^t \frac{\text{Li}_2(t)}{1-t} dt.$$

Clearly,

$$\int_0^t \frac{\text{Li}_2(t)}{t} dt = \text{Li}_3(t),$$

and integrating by parts,

$$\int_0^t \frac{\text{Li}_2(t)}{1-t} dt = -\ln(1-t)\text{Li}_2(t) - \int_0^t \frac{\ln^2(1-t)}{t} dt,$$

so that

$$\sum_{n=1}^{\infty} H_n^{(2)} \frac{t^n}{n} = \text{Li}_3(t) - \ln(1-t)\text{Li}_2(t) - \int_0^t \frac{\ln^2(1-t)}{t} dt \qquad (5)$$

At the same time

$$\sum_{n=1}^{\infty} H_n t^n = \frac{-\ln(1-t)}{1-t}.$$

Dividing this by $t$ and integrating we find

$$\sum_{n=1}^{\infty} H_n \frac{t^n}{n} = \int_0^t \frac{-\ln(1-t)}{t(1-t)} dt = \int_0^t \frac{-\ln(1-t)}{1-t} dt + \int_0^t \frac{-\ln(1-t)}{t} dt,$$

that is,

$$\sum_{n=1}^{\infty} H_n \frac{t^n}{n} = \frac{1}{2}\ln^2(1-t) + \text{Li}_2(t).$$

Multiplying this by 2 and then dividing by $t$ we find

$$2\sum_{n=1}^{\infty} H_n \frac{t^{n-1}}{n} = \frac{\ln^2(1-t)}{t} + 2\sum_{n=1}^{\infty} \frac{t^{n-1}}{n^2},$$

which integrated becomes

$$2\sum_{n=1}^{\infty} H_n \frac{t^n}{n^2} = \int_0^t \frac{\ln^2(1-t)}{t} dt + 2\text{Li}_3(t). \qquad (6)$$

Adding equations (5) and (6) gives the desired result.

We continue now with the evaluation of the series

$$\sum_{n=1}^{\infty} \frac{(-1)^{n+1}}{n+1} H_n^{(2)} = \frac{\pi^2}{12} \ln 2 - \frac{1}{4} \zeta(3).$$

Setting $t = -1$ in (4) yields

$$\sum_{n=1}^{\infty} H_n^{(2)} \frac{(-1)^{n-1}}{n} + 2 \sum_{n=1}^{\infty} H_n \frac{(-1)^{n-1}}{n^2} = 3 \operatorname{Li}_3(-1) - \operatorname{Li}_2(-1) \ln 2 = \frac{9}{4} \zeta(3) - \frac{\pi^2}{12} \ln 2.$$

The following evaluation is due to Sitaramachandra Rao [7]

$$\sum_{n=1}^{\infty} H_n \frac{(-1)^{n-1}}{n^2} = \frac{5}{8} \zeta(3).$$

Therefore, the above equation implies

$$\sum_{n=1}^{\infty} H_n^{(2)} \frac{(-1)^{n-1}}{n} = \zeta(3) - \frac{\pi^2}{12} \log 2.$$

From here we easily compute

$$\sum_{n=1}^{\infty} \frac{(-1)^{n+1}}{n+1} H_n^{(2)} = -\sum_{n=1}^{\infty} \frac{(-1)^{n-1}}{n} H_n^{(2)} - \operatorname{Li}_3(-1) = \frac{\pi^2}{12} \ln 2 - \frac{1}{4} \zeta(3)$$

which is the value of $J_2$ needed in section 2.